\documentclass{article}

\usepackage{amssymb,amsmath,theorem,euscript}
\input{epsf}

\newcounter{sec}

\newcounter{punct}[sec]

\def\punct{\refstepcounter{punct}{\arabic{sec}.\arabic{punct}.  }\boldmath}

\def\COUNTERS{\addtocounter{sec}{1}
              \setcounter{punct}{0}
          \setcounter{equation}{0}
          \setcounter{theorem}{0}
                  }

\newtheorem{theorem}{Theorem}[sec]
\newtheorem{proposition}[theorem]{Proposition}

\newtheorem{lemma}[theorem]{Lemma}

\newtheorem{corollary}[theorem]{Corollary}

 \def\ov{\overline}
\def\wt{\widetilde}
\def\wh{\widehat}

     \newcommand{\URep}{\mathop {\mathrm {URep}}\nolimits}
     \newcommand{\Rep}{\mathop {\mathrm {Rep}}\nolimits}
     \newcommand{\cRes}{\mathop {\mathrm {\cR es}}\nolimits}
 
\begin{document}

\def\OO{\mathrm{O}}
\def\GLO{\mathrm{GLO}}
\def\Coll{\mathrm{Coll}}
\def\kappa{\varkappa}
\def\Mat{\mathrm{Mat}}
\def\U{\mathrm U}
\def\GLL{\overline{\mathrm {GL}}}
\def\Spp{\overline{\mathrm {Sp}}}
\def\SL{\mathrm{SL}}
\def\SO{\mathrm{SO}}
\def\PGL{\mathrm{PGL}}

\def\R{\mathbb{R}}
\def\C{\mathbb{C}}
\def\V{\mathbb{V}}

\def\la{\langle}
\def\ra{\rangle}

 \def\cA{\mathcal A}
\def\cB{\mathcal B}
\def\cC{\mathcal C}
\def\cD{\mathcal D}
\def\cE{\mathcal E}
\def\cF{\mathcal F}
\def\cG{\mathcal G}
\def\cH{\mathcal H}
\def\cJ{\mathcal J}
\def\cI{\mathcal I}
\def\cK{\mathcal K}
 \def\cL{\mathcal L}
\def\cM{\mathcal M}
\def\cN{\mathcal N}
 \def\cO{\mathcal O}
\def\cP{\mathcal P}
\def\cQ{\mathcal Q}
\def\cR{\mathcal R}
\def\cS{\mathcal S}
\def\cT{\mathcal T}
\def\cU{\mathcal U}
\def\cV{\mathcal V}
 \def\cW{\mathcal W}
\def\cX{\mathcal X}
 \def\cY{\mathcal Y}
 \def\cZ{\mathcal Z}
 
 \def\cGL{\mathcal{GL}}
 
 \def\EXP{\mathcal{EXP}}
 
%%% END MATHCAL %%%%%%%%%%%%%%%%%%%%%%%%%%%%%%%%% %%%%%%%%%%%%%%%%%%%%%%%%%%%%%%%% %%%
%\def\0{{\ov 0}}
% \def\1{{\ov 1}}
 %%%%%%%%%%%%%%%%%%%%%%%%%%%% %%%%%%%%%%%%%%%%%%%%%%%%%%%%%%%%%%% %%% BEGIN GOTIC
 \def\frA{\mathfrak A}
 \def\frB{\mathfrak B}
\def\frC{\mathfrak C}
\def\frD{\mathfrak D}
\def\frE{\mathfrak E}
\def\frF{\mathfrak F}
\def\frG{\mathfrak G}
\def\frH{\mathfrak H}
\def\frI{\mathfrak I}
 \def\frJ{\mathfrak J}
 \def\frK{\mathfrak K}
 \def\frL{\mathfrak L}
\def\frM{\mathfrak M}
 \def\frN{\mathfrak N} \def\frO{\mathfrak O} \def\frP{\mathfrak P} \def\frQ{\mathfrak Q} \def\frR{\mathfrak R}
 \def\frS{\mathfrak S} \def\frT{\mathfrak T} \def\frU{\mathfrak U} \def\frV{\mathfrak V} \def\frW{\mathfrak W}
 \def\frX{\mathfrak X} \def\frY{\mathfrak Y} \def\frZ{\mathfrak Z} \def\fra{\mathfrak a} \def\frb{\mathfrak b}
 \def\frc{\mathfrak c} \def\frd{\mathfrak d} \def\fre{\mathfrak e} \def\frf{\mathfrak f} \def\frg{\mathfrak g}
 \def\frh{\mathfrak h} \def\fri{\mathfrak i} \def\frj{\mathfrak j} \def\frk{\mathfrak k} \def\frl{\mathfrak l}
 \def\frm{\mathfrak m} \def\frn{\mathfrak n} \def\fro{\mathfrak o} \def\frp{\mathfrak p} \def\frq{\mathfrak q}
 \def\frr{\mathfrak r} \def\frs{\mathfrak s} \def\frt{\mathfrak t} \def\fru{\mathfrak u} \def\frv{\mathfrak v}
 \def\frw{\mathfrak w} \def\frx{\mathfrak x} \def\fry{\mathfrak y} \def\frz{\mathfrak z} \def\frsp{\mathfrak{sp}}
 %% This is Lie algebra %%% END GOTIC
%%%%%%%%%%%%%%%%%%%%%%%%%%%%%%%% %%%%%%%%%%%%%%%%%%%%%%%%%%%%%%%%%
%%% BEGIN MATHBF
 \def\bfa{\mathbf a} \def\bfb{\mathbf b} \def\bfc{\mathbf c} \def\bfd{\mathbf d} \def\bfe{\mathbf e} \def\bff{\mathbf f}
 \def\bfg{\mathbf g} \def\bfh{\mathbf h} \def\bfi{\mathbf i} \def\bfj{\mathbf j} \def\bfk{\mathbf k} \def\bfl{\mathbf l}
 \def\bfm{\mathbf m} \def\bfn{\mathbf n} \def\bfo{\mathbf o} \def\bfp{\mathbf p} \def\bfq{\mathbf q} \def\bfr{\mathbf r}
 \def\bfs{\mathbf s} \def\bft{\mathbf t} \def\bfu{\mathbf u} \def\bfv{\mathbf v} \def\bfw{\mathbf w} \def\bfx{\mathbf x}
 \def\bfy{\mathbf y} \def\bfz{\mathbf z} \def\bfA{\mathbf A} \def\bfB{\mathbf B} \def\bfC{\mathbf C} \def\bfD{\mathbf D}
 \def\bfE{\mathbf E} \def\bfF{\mathbf F} \def\bfG{\mathbf G} \def\bfH{\mathbf H} \def\bfI{\mathbf I} \def\bfJ{\mathbf J}
 \def\bfK{\mathbf K} \def\bfL{\mathbf L} \def\bfM{\mathbf M} \def\bfN{\mathbf N} \def\bfO{\mathbf O} \def\bfP{\mathbf P}
 \def\bfQ{\mathbf Q} \def\bfR{\mathbf R} \def\bfS{\mathbf S} \def\bfT{\mathbf T} \def\bfU{\mathbf U} \def\bfV{\mathbf V}
 \def\bfW{\mathbf W} \def\bfX{\mathbf X} \def\bfY{\mathbf Y} \def\bfZ{\mathbf Z} \def\bfw{\mathbf w}
 %%% END MATHBF
%%%%%%%%%%%%%%%%%%%%%%%%%%%%%%% %%%%%%%%%%%%%%%%%%%%%%%%%%%%%%%%%
\def\B{\mathrm{B}}
%%%%%%%%%%%%%%%%%%%%%%%%%%%%%%%
 %%% BEGIN MATHBB
 \def\R {{\mathbb R }} \def\C {{\mathbb C }} \def\Z{{\mathbb Z}} \def\H{{\mathbb H}}
  \def\K{{\mathbb K}}
   \def\k{{\Bbbk}}
 \def\N{{\mathbb N}} \def\Q{{\mathbb Q}} \def\A{{\mathbb A}} \def\T{\mathbb T} 
 \def\G{\mathbb G}
 \def\bbA{\mathbb A} \def\bbB{\mathbb B} \def\bbD{\mathbb D} \def\bbE{\mathbb E} \def\bbF{\mathbb F} \def\bbG{\mathbb G}
 \def\bbI{\mathbb I} \def\bbJ{\mathbb J} \def\bbL{\mathbb L} \def\bbM{\mathbb M} \def\bbN{\mathbb N} \def\bbO{\mathbb O}
 \def\bbP{\mathbb P} \def\bbQ{\mathbb Q} \def\bbS{\mathbb S} \def\bbT{\mathbb T} \def\bbU{\mathbb U} \def\bbV{\mathbb V}
 \def\bbW{\mathbb W} \def\bbX{\mathbb X} \def\bbY{\mathbb Y} \def\kappa{\varkappa} \def\epsilon{\varepsilon}
 \def\phi{\varphi} \def\le{\leqslant} \def\ge{\geqslant}

\def\P{\mathbf P}

\def\GL{\mathrm {GL}}
\def\bGL{\mathbf {GL}}
\def\GLB{\mathrm {GLB}}

\def\bGr{\mathbf {Gr}}
\def\Gr{\mathrm {Gr}}
\def\Sp{\mathrm {Sp}}
\def\bFl{\mathbf {Fl}}

\def\1{\mathbf {1}}
\def\0{\mathbf {0}}

\def\rra{\rightrightarrows}

 \newcommand{\Dim}{\mathop {\mathrm {Dim}}\nolimits}
  \newcommand{\codim}{\mathop {\mathrm {codim}}\nolimits}
   \newcommand{\im}{\mathop {\mathrm {im}}\nolimits}
\newcommand{\ind}{\mathop {\mathrm {ind}}\nolimits}
\newcommand{\graph}{\mathop {\mathrm {graph}}\nolimits}

\def\F{\bbF}

\def\lambdA{{\boldsymbol{\lambda}}}
\def\alphA{{\boldsymbol{\alpha}}}
\def\betA{{\boldsymbol{\beta}}}
\def\gammA{{\boldsymbol{\gamma}}}
\def\deltA{{\boldsymbol{\delta}}}
\def\mU{{\boldsymbol{\mu}}}
\def\nU{{\boldsymbol{\nu}}}
\def\epsiloN{{\boldsymbol{\varepsilon}}}
\def\phI{{\boldsymbol{\phi}}}
\def\psI{{\boldsymbol{\psi}}}
\def\kappA{{\boldsymbol{\kappa}}}

\def\sm{\smallskip}
\def\nw{\nwarrow}
\def\se{\searrow}

\def\kos{/\!\!/}

\begin{center}
\Large\bf

Some remarks on traces on the infinite-dimensional Iwahori--Hecke algebra 

\medskip

\large \sc
Yury A. Neretin%
\footnote{Supported by the grant of FWF (Austrian science fund) P31591.}

\end{center}

{\small 
The infinite-dimensional Iwahori--Hecke algebras $\mathcal{H}_\infty(q)$ are direct
limits of the usual finite-dimensional Iwahori--Hecke algebras.
They arise in a natural way as  convolution
algebras of bi-invariant functions on  groups $\mathrm{GLB}(\mathbb{F}_q)$
 of infinite-dimensional
matrices over finite-fields having only finite number of non-zero matrix
elements under the diagonal. In 1988 Vershik and Kerov classified all
indecomposable positive traces on $\mathcal{H}_\infty(q)$. Any such trace
generates a representation of the double $\mathcal{H}_\infty(q)\otimes \mathcal{H}_\infty(q)$
and of the double $\mathrm{GLB}(\mathbb{F}_q)\times \mathrm{GLB}(\mathbb{F}_q)$.
We present  constructions of such representations; the traces
are some distinguished matrix elements. We also obtain some (simple)
general statements on relations between unitary representations of groups 
and representations of convolution algebras of measures bi-invariant with respect to compact subgroups.}

\section{Introduction}

\COUNTERS 

%In \cite{VK0}, \cite{Ker}, \cite{VK1}, \cite{VK2}  Vershik and Kerov classified all
%indecomposable traces on the infinite-dimensional Iwahori--Hecke algebras $\cH_\infty(q)$.
%They also introduced
%the group $\GLB(\F_q)$ of infinite matrices over a finite field 
%$\F_q$ having only finite number of non-zero elements under the diagonal, and
% showed that representations of the algebra $\cH_\infty(q)$ correspond to type $\mathrm{II_\infty}$
%or $\mathrm{I_\infty}$
%factor-representations of the group $\GLB(\F_q)$, \cite{VK1}. Because of
%the death of Kerov, the work was not completely published.
% See also continuations of this work
%in \cite{Mel}, \cite{BP},  \cite{GKV}.  The present paper contains an explicit
%construction of  representations of doubles $\cH_\infty(q)\otimes\cH_\infty(q)$
%corresponding to indecomposable traces of $\cH_\infty(q)$, and some formal
%logical arrows, which can be unclear from previous publications
% (in particular, we present quasi-explicit constructions
%of unitary representations of $\GLB(\F_q)\times \GLB(\F_q)$ corresponding to the indecomposable traces).

%\sm 

{\bf\punct Algebras of bi-invariant measures.%
\label{ss:biinvariant}}
Let $G$ be a topological group. Denote by $\cM(G)$ the algebra of finite compactly supported
complex-valued measures on $G$. The addition in $\cM(G)$ is the  addition of measures, the multiplication
is the convolution. Namely, let $\mu_1$, $\mu_2$ be measures supported by compact sets $L_1$, $L_2$ respectively. Consider the measure $\mu_1\times \mu_2$ on $L_1\times L_2$. The {\it convolution} $\mu_1*\mu_2$
is the pushforward of $\mu_1\times\mu_2$ under the map $(g_1,g_2)\mapsto g_1g_2$ from $L_1\times L_2$
to $G$. We also define an {\it involution} $\mu\mapsto \mu^\star$ in $\cM(G)$. Namely,
$\mu^\star$ is the pushforward of the the complex conjugate measure $\mu$ under the map
$g\mapsto g^{-1}$,
$$
(\mu_1*\mu_2)^\star=\mu_2^\star *\mu_1^\star.
$$

We also define a {\it transposition} that send $\mu$ to its image under the map
$g\mapsto g^{-1}$.  We have
$$
(\mu_1*\mu_2)^t=\mu_2^t *\mu_1^t,
$$
the involution is anti-linear map and the transposition is linear.

\sm 

{\sc Remark.} If a group $G$ is finite, then $\cM(G)$ is the group algebra
of $G$.
\hfill $\boxtimes$

\sm 

 Let $\rho$ be a unitary representation of $G$ in a Hilbert space
$H$ ({\it Hilbert spaces assumed to be separable}). For any measure $\mu\in \cM(G)$ denote by
$\rho(\mu)$ the operator 
$$
\rho(\mu)=\int_G \rho(g)\,d\mu(g).
$$
This determines a {\it $*$-representation} of the algebra $\cM(G)$,
$$
\rho(\mu_1)\rho(\mu_2)=\rho(\mu_1*\mu_2), \qquad \rho(\mu)+\rho(\nu)=\rho(\mu+\nu),
\qquad \rho(\mu^\star)=\rho(\mu)^*.
$$

Let $K\subset G$ be a compact subgroup. By $\delta_K$ we denote the probabilistic Haar measure on
 $K$ regarded as an element of $\cM(G)$. Denote by $\cM(G\kos K)$ the subalgebra of $\cM(G)$ consisting
 of measures, which are invariant with respect to left and right shifts by elements of $K$.
 Clearly, for any $\mu\in \cM(G)$ we have $\delta_K*\mu*\delta_K\in \cM(G\kos K)$
 and for an element $\nu\in \cM(G)$ the following identity holds:
 $$
 \nu\in\cM(G\kos K)\Longleftrightarrow \delta_K*\nu*\delta_K=\nu.
 $$
 Clearly, the subalgebra $\cM(G\kos K)$ is closed with respect to
 the involution.
 
 \sm 
 
 Let $\rho$ be a unitary representation of $G$ in a Hilbert space $H$. Denote
 by $H^K\subset H$ the subspace consisting of all $K$-fixed vectors, denote by $P^K$
 the operator of orthogonal projection to $H^K$.
 Clearly,
$$
P^K:=\rho(\delta_K).
$$
Therefore, for any $\nu\in \cM(G)$ we have 
$$P^K\rho(\nu)P^K=\rho(\nu).$$
Therefore, such operators have the following block form
with respect to the decomposition $H=H^K\oplus (H^K)^\bot$:
\begin{equation}
\rho(\nu)=
\begin{pmatrix}
\wt \rho(\nu)&0\\0&0
\end{pmatrix}.
\label{eq:000}
\end{equation}
So for any unitary representation $\rho$ of $G$
 we get  a $*$-representation $\wt\rho(\cdot)$ of the algebra $\cM(G\kos K)$
in the Hilbert space $H^K$.

Denote by $\URep(G)_K$ the set of equivalence classes of unitary representations of the group $G$ such that $H^K$ is a cyclic subspace%
\footnote{I.e., linear combinations of vectors $\rho(g)\xi$, where $\xi$ ranges in $H^K$ and
	$g$ ranges in $G$, are dense in $H$.}.
 Denote by $\Rep(G\kos K)$ the set of $*$-representations
	of the algebra $\cM(G\kos K)$. The following statement is obvious
	(see below Proposition \ref{pr:1})
	
	\sm 
	
{\it The map $\cRes_{G\kos K}:\rho\mapsto \wt\rho$ from $\URep(G)_K$ to $\Rep(\cM(G\kos K))$ is injective.}
			
	\sm

An inverse construction $\wt \rho\mapsto\rho$ (if  $\wt \rho$ is contained in the image of $\cR_{G\kos K}$) is semi-explicit: having a representation of $\cM(G\kos K)$
one can define a reproducing kernel space and a representation of $G$ in this space
(see Proposition \ref{pr:2}).
Generally speaking, a positive-definiteness of a kernel is a non-trivial question  and
a description of the image of $\cRes_{G\kos K}$ also is non-trivial%
\footnote{For instance, if $G$ is a semisimple Lie group and $K$ is its maximal compact subgroup.}.
For our purposes the following statement is sufficient:

\sm 

{\it If $G$ is a compact group or a direct limit of compact groups,
	 then the map $\cRes_{G\kos K}$ is a bijection},
 see Propositions \ref{pr:3}--\ref{pr:4}. 

\sm

\sm

Recall also two variations of these definitions.
Let $G$ be a unimodular locally compact group%
\footnote{A locally compact group is called {\it unimodular} if its Haar measure is two-side
invariant, see, e.g., \cite{Kir}, Subsect. 9.1.}, let $\lambda$ be a Haar measure on $G$. 
Let $K\subset G$ be a compact subgroup. Denote
by $\cC(G)$ the convolution algebra of all compactly supported continuous functions on $G$.
This algebra is a subalgebra in $\cM(G)$, namely, for any function $\phi\in \cC(G)$ we assign
the measure $\phi(g)\,\mu(g)$. By $\cC(G\kos K)\subset \cC(G)$ we denote the subalgebra
of $K$-bi-invariant functions, $\phi(k_1 gk_2)=\phi(g)$ for $k_1$, $k_2\in K$.

Next, let $G$ be a unimodular totally disconnected locally compact group, $k$ be an open compact subgroup.
Denote by $\cA(G\kos K)$ the convolution algebra consisting of $K$-bi-invariant compactly
supported locally constant functions. Clearly, in this case
$$\cA(G\kos K)= \cC(G\kos K)=\cM(G\kos K). 
$$

Convolution algebras  $\cM(G\kos K)$ are a usual tool of representation theory,
see, e.g., \cite{Kre}, \cite{Gel}, \cite{BG}, \cite{Iwa}, \cite{IwaM}, \cite{CIK}, \cite{BZ}, \cite{Koo}, \cite{Macd}, Chapter 5. 

\sm 

{\bf \punct Iwahori-Hecke algebras.%
\label{ss:iwahori}}
Consider a finite field $\F_q$ with $q$ elements.
Let $G$ be the group $\GL(n,\F_q)$ of invertible matrices of order $n$ over $\F_q$.
Let $K=B(n)$ be the subgroup of upper triangular matrices. The algebra $\cM(G\kos K)$ was described by
Iwahori \cite{Iwa}. Denote by $s_m\in \GL(n,\F_q)$ the operator that transposes basis vectors $e_m$ and $e_{m+1}$
in $\F_q^n$ and fixes other $e_j$. So $m=1$, 2, \dots, $n-1$. Set $\sigma_m:=\delta_K*s_m*\delta_K$. Then 
the elements $\sigma_m$
 generate the algebra 
$$
\cH_n(q):=\cM\bigl(\GL(n,\F_q)\kos B(n)\bigr)
$$
and relations are 
\begin{align}
\sigma_m^2&=(q-1)\sigma_k+q\qquad\text{or}\qquad (\sigma_m+1)(\sigma_m-q)=0;
\label{eq:h-1}
\\
\sigma_m \sigma_{m+1} \sigma_m&=\sigma_{m+1} \sigma_m\sigma_{m+1}
\label{eq:h-2}
\\
\sigma_m\sigma_l&=\sigma_l\sigma_m \qquad \text{if $|m-l|>1$.}
\label{eq:h-3}
\end{align}
The dimension of $\cH_q=n!$. The involution and the transposition are determined by 
$$\sigma^\star_m=\sigma_m, \qquad \sigma^t_m=\sigma_m.$$

By the construction $q$ is a power of a prime. However, the algebra with
relations \eqref{eq:h-1}--\eqref{eq:h-3} makes sense for any $q\in\C$, for $q=1$
we get the group algebra of the symmetric group $S_n$. For all $q$ that are not roots of units
algebras $\cH_n(q)$ are isomorphic. Therefore they have the same dimensions of representations.

\sm 

{\bf\punct The group of almost triangular matrices and its Iwahori--Hecke algebra.%
\label{ss:iwahori-infinity}}  
Following Vershik and Kerov \cite{VK1}, denote by $\GLB(\F_q)$ the group of all infinite invertible matrices over $\F_q$ having
only finite number of nonzero matrix elements under the diagonal%
\footnote{There are several approaches to 'representation theory of infinite-dimensional
groups over finite fields', see a discussion of different  works  in \cite{Ner-finite}, Subsect. 1.11.}. 

Denote by $B(\infty)$
the subgroup of $\GLB(\F_q)$ consisting of upper triangular matrices. 
For a matrix $g\in K$ its diagonal elements $g_{jj}$ are contained in $\F_q^\times:=\F_q\setminus 0$
and elements $g_{ij}$, where $i<j$ are contained in $\F_q$. So the set $B(\infty)$ is a product
of a countable number of copies if $\F_q^\times$ and a countable number of copies of $\F_q$.
We equip $K$ with the product topology and get a structure of a compact topological group.
We take uniform probabilistic measures on the sets $\F_q^\times$, $\F_q$
and equip the group $B(\infty)$ with the product measure.  It is easy to verify that we get 
the Haar measure on $B(\infty)$.  

The space 
 $\GLB(\F_q)$ is a disjoint union of a countable number of  cosets $g B(\infty)$. We equip $\GLB(\F_q)$ with the topology
 of disjoint union, this determines 
a structure of a unimodular locally compact topological group on $\GLB(\F_q)$.
Equivalently, a sequence $g^{(\alpha)}\in\GLB(\F)$ converges to $g$
if the following two conditions hold: 

\sm

--- for each $i$, $j$ we have $g_{ij}^\alpha=g_{ij}$ for sufficiently large $\alpha$;

\sm

--- there exists $n$ such that for all $i$, $j$ such that $i<j$ and $j\ge n$ we have
$g_{ij}^\alpha=0$ for all $\alpha$. 

\sm

Denote by $\GLB(n,\F_q)\subset \GLB(\F_q)$ the subgroup consisting of matrices $g$ satisfying the condition:
$g_{ij}=0$ for all pairs $(i,j)$ such that $j>n$, $i<j$. In other words this group
is generated by $\GL(n,\F_q)$ and $B(\infty)$.
 We get an increasing family
of compact subgroups in $\GLB(\F_q)$, and 
$\GLB(\F_q)$ is the direct limit
$$
\GLB(\F_q)=\lim_{n\longrightarrow\infty} \GLB(n,\F_q).
$$

\sm 

The Iwahori--Hecke algebra 
$$\cH_\infty(q):=\cM\bigl(\GLB(\F_q)\kos B(\infty_q)\bigr)=\cA\bigl(\GLB(\F_q)\kos B(\infty_q)\bigr)$$
 is the algebra with generators $\sigma_1$, $\sigma_2$, $\sigma_3$, \dots
and the same relations \eqref{eq:h-1}--\eqref{eq:h-3}. Also, it is the direct limit
$$
\cH_\infty(q)=\lim_{n\longrightarrow \infty } \cH_n(q).
$$

Again, algebras $\cH_\infty(q)$ are well defined for all $q\in\C$.
By default we assume
$$
q>0.
$$

\sm

{\bf \punct Traces on $\cH_\infty(q)$.%
\label{ss:traces}}
 A {\it trace} $\chi$ on an associative algebra $\cA$ with involution is a linear functional
$\cA\to \C$ such that

\sm 

a) $\chi(A B)=\chi(B A)$ for any $A$, $B\in \cA$;

\sm 

b) $\chi(A^\star)=\ov {\chi(A)}$.

\sm

c) $\chi(A^\star A)\ge 0$ for any $A\in \cA(q)$.

\sm 

The set of all traces is a convex cone: if $\chi_1$, $\chi_2$ are traces, then
for any $a_1$, $a_2\ge 0$ the expression $a_1\chi_1+a_2\chi_2$ is a trace.
A trace $\chi$ is {\it normalized} if $\chi(1)=1$. 
A trace $\nu$ is {\it indecomposable} if it does not admit a representation
of the form $\nu= a_1\chi_1+a_2\chi_2$, where $\chi_1$, $\chi_2$ are traces
non-proportional to $\nu$ and $a_1$, $a_2>0$.

\sm 

In \cite{VK0}, \cite{Ker}, \cite{VK1}, \cite{VK2}  Vershik and Kerov classified%
\footnote{Kerov died in 2020, the work was not completely published, the text \cite{VK2}
	was based on his posthumos notes.} all
indecomposable traces on the infinite-dimensional Iwahori--Hecke algebras $\cH_\infty(q)$.
See also continuations of this work
in \cite{Mel}, \cite{BP},  \cite{GKV}.

Let us formulate the classification theorem.
For $\lambda\le \nu$ we denote 
$$
\zeta_{\lambda,\mu}:=\sigma_{\mu-1}\sigma_{\mu-2}\dots \sigma_\lambda.
$$
Consider a partition $\nU$ of $n$,
$$
n=\nu_1+\nu_2+\dots , \qquad \nu_1\ge\nu_2\ge\dots >0.
$$
We set
$$
\lambda_j:=\nu_1+\dots +\nu_j,
$$
and consider the following elements of the Iwahori--Hecke algebra
$$
\zeta_{\lambdA}:= \dots \zeta_{\lambda_3,\lambda_4} \, \zeta_{\lambda_2,\lambda_3}   \, \zeta_{\lambda_1,\lambda_2}
$$
(all factors commute). For the following two statements, see \cite{Ram}, \cite{GP}, Sect.8.2.

\sm

--- {\it The elements $\zeta_\lambdA$ form a basis of the space}
$$
\cH_n(q)/[\cH_n(q),\cH_n(q)],
$$ 
where   $[\cH_n(q),\cH_n(q)]\subset \cH_n(q)$ is the subspace generated by all commutators $ab-ba$,
where $a$, $b$ range in $\cH_n(q)$. Therefore,

\sm 

--- {\it a trace on $\cH_n(q)$ is determined by its values on elements $\zeta_{\lambdA}$.}

\sm 

This implies that the same statement is valid for $\cH_\infty(q)$.  Next,

\sm

---  By \cite{VK0}, {\it  for any indecomposable trace $\chi$ on $\cH_\infty(q)$
we have
\begin{equation}
\label{eq:multi}
\chi(\zeta_{\lambdA})=\prod_j \zeta_{\lambda_j},
\end{equation}
and for any $m>0$}
\begin{equation}
\chi(\zeta_{[\lambda_1+m,\lambda_2+m]})=\chi(\zeta_{[\lambda_1,\lambda_2]})
.
\label{eq:shift}
\end{equation}

Therefore, any trace $\chi$ on $\cH_\infty(q)$ is determined
by its values on elements
$$
\zeta_m:=\zeta_{[1,m]}.
$$
Notice, that $\zeta_1=1$.

% Vershik and Kerov 
%\cite{VK0}--\cite{VK1} classified all indecomposable traces on the algebra
%$\cH_\infty(q)$ with $q>0$.

\begin{theorem}
	\label{th:vk}
	 {\rm(Vershik, Kerov \cite{VK0})} Let $q>0$.
 Indecomposable traces $\chi^{\alpha,\beta,\gamma}$ on $\cH_\infty(q)$ are enumerated by
collections of parameters
\begin{align*}
&\alpha_1\ge \alpha_2\ge \dots \ge 0,
\qquad 
\beta_1\ge \beta_2\ge \dots\ge 0,\qquad \gamma\ge0;
\\
&\sum \alpha_i +\sum \beta_j+\gamma=1.
\end{align*}
The value of $\chi^{\alpha,\beta,\gamma}$ on $\zeta_m$
is given by the formula
\begin{equation}
\chi^{\alpha,\beta,\gamma}(\zeta_m)=
\frac1{q-1}
\sum_{\text{$\begin{matrix}\mu_1\ge 0,\mu_2\ge 0,\dots 
	\\
\sum k\mu_k=m
\end{matrix}$}
}  \,\, \prod_{k\ge 1} \frac{ (q^k-1)^{\mu_k}}
	{k^{\mu_k} \mu_k! } \prod_{k\ge2} p_k(\alpha,\beta)^{\mu_k},
	\label{eq:vk}
\end{equation}
where
\begin{equation}
p_k(\alpha,\beta):= \sum_i \alpha_i^k+(-1)^{k+1}\sum_i\beta_i^k
\end{equation}
denote the {\rm super-Newton sums}.
	\end{theorem}

{\sc Remark.} Recall that symmetric functions \cite{Macd} can be represented as polynomials
of the Newton sums $p_k(\alpha):=\sum \alpha_i^k$. {\it Supersymmetric functions} are polynomials
in super-Newton sums $p_k(\alpha,\beta)$, see \cite{KOO}, \cite{Ser}, \cite{BO}, Chapter 2.

\sm 
 
{\sc Remark.}
For $q=1$ (i.e., for the infinite symmetric group)
the expression \eqref{eq:vk} has a removable singularity,
$$
\chi^{\alpha,\beta,\gamma}(\zeta_m)=p_m(\alpha,\beta)\qquad\text {for $m\ge 2$}.
$$ 
This special case of Theorem {\rm \ref{th:vk}} is the Thoma theorem \cite{Tho-symm} (see, also \cite{BO}, Chapter 4): irreducible normalized characters of the infinite symmetric group
 have the following
form
$$
\chi_{\alpha,\beta,\gamma}(g)=\prod_{m\ge 2} p_m(\alpha,\beta)^{r_m(g)},
$$ 
where $r_m(g)$ is the number of cycles of $g$ of order $m$.
\hfill $\boxtimes$

\sm 

{\bf\punct Traces and representations of the double $\cH_{\infty}(q)\otimes \cH_{\infty}(q)$.%
\label{ss:double}}
Let $\cA$ be a $*$-algebra with unit, let  $\chi$ be a trace on $\cA$.
Denote by $\cA^\circ$ the algebra anti-isomorphic to
$\cA$, it  coincides with $\cA$ as a linear space,
multiplication is given by $A\diamond B:=BA$. 
The following formula determines an inner product
on $\cA$:
$$
\la A,B\ra_\chi= \chi(B^\star A).
$$ 
Denote by $\ov \cA_\chi$ the corresponding Hilbert space.
Assume that  operators of left multiplication are
bounded%
\footnote{So have
a structure of a Hilbert algebra in the sense of Dixmier (see \cite{Dix-Hilbert}, \S I.6, \cite{Dix},  A.54).} in $\ov A_\chi$.
Then we have an action of $\cA$ on $\ov \cA$ by left multiplications and the action
 of $\cA^\circ$ by right multiplications. These actions commute, so we get a representation
 of the tensor product $\cA\otimes \cA^\circ$,
 $$
 \tau_\chi(A\otimes B)X:= AXB ,\qquad \text{where $X\in \ov\cA$, $A\in\cA$, $B\in\cA^\circ$.}
 $$

{\sc Remark.} Representations $\tau_\chi$ of $\cA\otimes \cA^\circ$ are not arbitrary representations,
they satisfy the condition:
$$
\tau_\chi(A\otimes 1)\, 1=\tau_\chi(1\otimes A)\, 1. \qquad \boxtimes
$$

\sm 

Notice, that 
\begin{equation}
\chi(A)= \bigl\la\tau_\chi(A\otimes 1) \,1,1\bigr\ra_\chi
\label{eq:1-1}
.\end{equation}
So a trace on $\cA$ is a certain matrix element
of a certain  $(\cA,\cA^\circ)$-bimodule.

\sm

Let us return to a discussion the Iwahori-Hecke algebras.
Notice that a transposition in $\cH_{\infty}(q)$
is an anti-isomorphism, so we can regard $\ov{\cH_{\infty}(q)}$ as a
$(\cH_{\infty}(q)\otimes\cH_{\infty}(q))$-module.

\begin{lemma} We have an isomorphism of algebras
	\begin{multline*}
	\cM\Bigl(\GLB(\F_q)\times \GLB(\F_q)\kos B(\infty)\times B(\infty)\Bigr) 
	\simeq\\
	\simeq \cM\Bigl(\GLB(\F_q)\kos B(\infty)\Bigr)\otimes 
	\cM\Bigl(\GLB(\F_q)\kos B(\infty)\Bigr)=\cH_\infty(q)\otimes \cH_\infty(q).
	\end{multline*}
	
	More generally, for any locally compact group $G$ with a countable base of topology
	and an open compact subgroup $K$, the following algebras are isomorphic
	$$
	\cM\bigl(G\times G\kos K\times K \bigr)\simeq \cM(G\kos K)\otimes \cM(G\kos K).
	$$
\end{lemma}

\begin{lemma}
	\label{l:bounded}
	For any indecomposable trace $\chi^{\alpha,\beta,\gamma}$ on $\cH_{\infty}(q)$
	 the operators of left-right multiplication
	  $$
	 \tau_\chi(A\otimes B)\,X:= AXB
	 $$
	are bounded in $\ov{\cH_{\infty}(q)}_\chi$.
\end{lemma}

The first statement is obvious,
the second statement is clear from the condition
\eqref{eq:h-1}. Indeed, the spectrum of a self-adjoint operator
of multiplication by a generator $\sigma_m$ consists of points $1$ and $q$, so its norm
is $\max(1,q)$. 

Our Proposition \ref{pr:4} implies the following corollary:

\begin{corollary}
	\label{cor:double}
	{\rm(}see, {\rm\cite{VK1})}
	Let $q=p^l$, where $p$ is a prime. Then for any indecomposable trace $\chi^{\alpha,\beta,\gamma}$ there
	exists an irreducible representation $\rho$ of the double $\GLB(\F_q)\times \GLB(\F_q)$,
	for which
	$\cRes(\rho)$ is isomorphic to the representation $\tau_\chi$ of $\cH_\infty(q)\otimes \cH_\infty(q)$.
\end{corollary}

\sm 

{\sc Remark.}
To avoid a misleading, we must say some  remarks.
For a type I topological group $G$ any irreducible representation of $G\times G$
is a tensor product $\pi_1\otimes \pi_2$ of irreducible representations of two copies of $G$
(see, e.g., \cite{Dix}, 13.1.8). The group $\GLB(\F_q)$ is not of type I and  a similar implication
is false. 
According \cite{VK0}, for an indecomposable trace $\chi$,
the operator algebra generated by
the representation of $\cH_{\infty}(q)\otimes 1$ in $V=\ov{\cH_{\infty}(q)}$
is a finite Murray--von Neumann factor. The representation of $1\otimes \cH_{\infty}(q)$
also generates a factor, which is the commutant of the first factor. The representation of the double
$\cH_{\infty}(q)\otimes \cH_{\infty}(q)$ is irreducible. The trace of the unit operator
is finite ($=1$), so it is a factor of a type $\mathrm {I}_n$ or $\mathrm {II}_1$.
But $\cH_{\infty}(q)$ has only two irreducible finite dimensional representations,
both are one-dimensional%
\footnote{Indeed, the algebra $\cH_n(q)$ is isomorphic to the group algebra of the symmetric group $S_n$,
so dimensions of their irreducible representations coincide. So they have two one-dimensional
representations, all other representations have dimensions $\ge n-1$. Therefore finite-dimensional 
irreducible representations of $\cH_\infty(q)$ are one-dimensional.}.
The first is the  trivial representation ($\alpha_1=1$, other $\alpha$'s and $\beta$'s are zero), the second is the representation
sending all $\sigma_j$ to  $-1$ ($\beta_1=1$, other $\beta$'s and $\alpha$'s are 0).  In the remaining cases we have $\mathrm {II}_1$-factors.
The  representations of $\GLB(\F_q)$ corresponding to the one-dimensional
characters are  the trivial representation 
and the Steinberg representation (see, e.g., \cite{Ner-steinberg}). The remaining characters
$\chi^{\alpha,\beta,\gamma}$ correspond to representations of $\GLB(\F_q)$ generating
 $\mathrm {II}_\infty$-factors.
\hfill $\boxtimes$

\sm 

Our next purpose is to construct explicitly  representations $\tau_\chi$ of 
$\cH_\infty(q)\otimes\cH_\infty(q)$. {\it We consider only the case
$\gamma=0.$}

\sm 
{\bf\punct $R$-matrix.%
\label{ss:R}}
Let $V$ be a Hilbert space equipped with an orthonormal basis
$v_j$, where $j$ ranges in non-zero integers.
Consider $R$-matrix given by the formula%
\footnote{For positive  indices $i$  the $R$ is a standard $R$-matrix (see, e.g., \cite{Jim}), the first summand in  formula \eqref{eq:R}
gives a kind of superization. For $q=1$ the matrix $R$ produces a structure equivalent
to a supersymmetric (or $\Z_2$-graded) tensor product,
see \cite{Ols-symm}, \cite{BO}, Chapter 10.}
\begin{multline}
R:=-\sum_{i<0} e_{ii}\otimes f_{ii} +q \sum_{i>0} e_{ii}\otimes f_{ii}-
\\-
\sqrt q \sum_{i\ne j,i\ne0, j\ne0} e_{ji}\otimes f_{ij} +(q-1) \sum_{i< j,i\ne0, j\ne0}
e_{ij}\otimes f_{ij},
\label{eq:R}
\end{multline}
see Kerov \cite{Ker}. Namely,  $R$ determines an operator in $V\otimes V$,
 here $e_{\alpha\beta}$ denote matrix units in
 the first copy of $V$,
 $$
 e_{\alpha\beta} v_i=\begin{cases}
 0,&\qquad \text{if $i\ne \alpha$};\\
 v_\beta,& \qquad \text{if $i\ne \alpha$},
 \end{cases}
 $$
$f_{\alpha\beta}$ denote matrix units in the second copy of $V$.
Then 
$$R^2=(q-1)R+q 1\otimes 1.$$
Consider the space $V^{\otimes n}=V\otimes V\otimes\dots$.
Denote by $R_{i(i+1)}$ the operator $R$ acting on the $i$-th and $(i+1)$-th factors.
Then 
%\begin{align*}
%R_{i(i+1)}R_{(i-1)i} R_{i(i+1)}=R_{(i-1)i} R_{i(i+1)} R_{(i-1)i};
%\\
%R_{i(i+1)} R_{j(j+1)}=R_{j(j+1)} R_{i(i+1)}, \qquad \text{for $|i-j|>1$.}
%\end{align*}
 the map 
 \begin{equation}
 \sigma_j\mapsto R_{j(j+1)}
 \label{eq:sigma-R}
 \end{equation}
  determines a representation
of $\cH_n(q)$ in $V^{\otimes n}$.

\sm 

{\bf \punct  Representations of $\cH_\infty(q)\otimes \cH_\infty(q)$.%
\label{ss:result}} 
Here we imitate the construction of representations of a double
of an infinite symmetric group proposed by Wassermann \cite{Was} and Olshanski \cite{Ols-symm}
(see, also, \cite{BO}, Chapter 10).

Let $\gamma=0$.
 Let $W$ be a copy of the space $V$ defined in the previous subsection,
let $w_j\in W$ be the basis corresponding to $v_j$. Consider the tensor
product $V\otimes W$ and the following unit vector\footnote{Recall that the definition of infinite tensor products of Hilbert spaces requires a
	choosing of distinguished vectors, see, e.g., \cite{Gui}, Addendum A.}
$$
\xi:=
\sum_{j>0} \sqrt{\beta_j}\, v_{-j}\otimes w_{-j}+\sum_{j>0} \sqrt{\alpha_j}\, v_{j}\otimes w_{j}\,\in V\otimes W.
$$
Since $\gamma=0$, we have $\|\xi\|^2=\sum \alpha_j+\sum \beta_j=1$.
Consider the infinite tensor product of Hilbert spaces
$$
\cX:=(V\otimes W,\xi)\otimes (V\otimes W,\xi)\otimes \dots 
$$
Denote
$$
\Xi:=\xi\otimes \xi\otimes \dots\in\cX.
$$
For $j\in \N$ define an operator 
$$R_{j(j+1)}^{\mathrm{left}}:\cX\to\cX$$
in the following way. We represent
$\cX$ as
\begin{multline*}
\cX=(V\otimes W,\xi)^{\otimes \infty}=\\=
(V\otimes W,\xi)^{\otimes (j-1)}\otimes \Bigl[(V\otimes W,\xi)\otimes (V\otimes W,\xi)\Bigr]
\otimes (V\otimes W,\xi)^{\otimes (\infty-j-1)}.
\end{multline*}
The middle factor is
$$V\otimes W\otimes V \otimes W \simeq (V\otimes V\otimes W\otimes W)$$
(for finite tensor products distinguished vectors have no matter). Then
$R_{j(j+1)}^{\mathrm{left}}$ is a tensor product of the following operators:

\sm 

--- $1^{\otimes(j-1)}$ in the first factor;

\sm 

---  $R\otimes (1\otimes 1)$ in the middle factor $(V\otimes V)\otimes (W\otimes W)$;

\sm 

--- $1^{\otimes(\infty-j-1)}$ in the last factor.

\sm

We get a representation of $\cH_\infty(q)$ in $\cX$, denote it by
$$A\mapsto A^{(l)}.$$
 In the same way 
we define operators $R_{j(j+1)}^{\mathrm{right}}$
acting by twisted permutations of factors $W$ and a representation of the second copy
of the algebra $\cH_\infty(q)$.

\begin{theorem}
	\label{th:}
	{\rm a)} For any $A\in\cH_\infty(q)$,
	$$
	\la A^{(l)}\Xi,\Xi\ra_{\cX}=\chi^{\alpha,\beta,0}(A).
	$$
	
	{\rm b)}
	Moreover, the representation of $\cH_\infty(q)\otimes \cH_\infty(q)$ in the cyclic span
	of $\Xi$ is equivalent to the representation in $\ov{\cH_\infty(q)}_\chi$.
\end{theorem}

{\bf\punct The further structure of the paper.} In Section 2 we show that
a unitary representation $\rho\in \URep(G)_K$ is uniquely determined by the
corresponding representation $\wt \rho$ of the algebra   $\cM(G\kos K)$ (Proposition \ref{pr:1})
and describe a quasi-explicit way of realization of $\rho$  (Proposition \ref{pr:2})
in a reproducing kernel space\footnote{On reproducing kernel spaces, see, e.g., \cite{Ner-gauss},
Sect. 7.1 and Subsect. 7.5.15.} constructed by $\wt\rho$. 
Next, we show that for a direct limit $L$ of compact groups any $*$-representation
of $\cM(L\kos K)$ corresponds to a unitary representation of $L$ (Proposition \ref{pr:4}). 
All these statements are obvious or very simple, however I could not find a source for
formal references. 

In Section 3 we prove Theorem \ref{th:}, this proof does not depend on Section 2.

\sm 

{\bf \punct Some comments.} The group $\GLB(\F_q)$ is a locally compact group, whose properties are partially similar to
properties of 'infinite dimensional' groups. The construction of Subsect. \ref{ss:result} looks like
a relatively usual construction related to infinite dimensional groups. However, a behavior of double coset
algebras is unusual. 

As an example of the usual behavior, consider $G_\alpha(n):=\GL(\alpha+n, \R)$, $K_\alpha(n):=\OO(n)$.
For $\alpha=0$ the algebra $\cC\bigl(G_0(n)\kos K_0(n)\bigr)$ is commutative and according Gelfand \cite{Gel}
 this implies 
sphericity of the subgroup $\OO(n)$ in the group $\GL(n,\R)$. However, even the algebra $\cC\bigl(\SL(2,\R)\kos \SO(2)\bigr)$
is a nontrivial object, a multiplication is defined in terms of a  hypergeometric
kernel,  see, e.g., \cite{Koo}. For $\alpha>0$ algebras 
$\cC\bigl(G_\alpha(n)\kos K_\alpha(n)\bigr)$ at the present time seem completely incomprehensible.

We have a natural map of double coset spaces
$$
K_\alpha(n)\setminus G_\alpha(n)/K_\alpha(n)\to K_\alpha(n+1)\setminus G_\alpha(n+1)/K_\alpha(n+1).
$$ 
but $K_\alpha(n+1)$ is strictly larger than $K_\alpha(n)$ and this map does not induce a homomorphism
from $\cC\bigl(G_\alpha(n)\kos K_\alpha(n)\bigr)$ to $\cC\bigl(G_\alpha(n+1)\kos K_\alpha(n+1)\bigr)$.
However, a limit object exists (see \cite{Olsh-GB}, \cite{Ner-concentraion}) and it has a structure of a semigroup.
Namely, there is a natural multiplication on the space of double cosets
$$
\Gamma_\alpha(\infty)=K_\alpha(\infty)\setminus G_\alpha(\infty)/K_\alpha(\infty),
$$
it admits a reasonable description (see \cite{Ner-book}, Section IX.4), and 
the semigroup $\Gamma_\alpha(\infty)$ acts in subspaces of $K_\alpha(\infty)$-fixed vectors
of unitary representations of $G_\alpha(\infty)=\GL(\infty,\R)$. Moreover, this situation is typical for infinite
dimensional groups and allows to produce unconventional algebraic structures (see, e.g., \cite{Ner-symm}). 

\sm 

In the context of the present paper, we have  increasing groups
$G(n)$ and  constant compact subgroup $K$, for this reason we have 
embeddings $\cM(G(n)\kos K)\to \cM(G(n+1)\kos K)$ and  a direct limit $\cM(G(\infty)\kos K)$.
Moreover, we have  comprehensive prelimit algebras $\cM(G(n)\kos K)$. This situation is unusual
among infinite-dimensional groups whose representation theories were topic of considerations of 
mathematicians. Certainly, additionally there is a group $\GL$ of two-side-infinite almost triangular matrices 
over $\F_q$ and its symplectic, orthogonal, and 'unitary' subgroups (in all cases the subgroup $K$
consists of upper triangular matrices), the author does not see a way to extend this list.

\section{Algebras of bi-invariant measures. Generalities}

\COUNTERS

%In this section we discuss a reconstruction
%of a unitary representation $\rho$ of a  group $G$
%from a representation $\tau$ of the algebra
%$\cM(G\kos K)$.
%We present a trivial Proposition \ref{pr:2} about a realization of $\rho$
%in a reproducing kernel space (on reproducing kernel spaces, see, e.g., \cite{Ner-gauss},
%Sect. 7.1 and Subsect. 7.5.15) and semi-trivial Propositions \ref{pr:3}--\ref{pr:4}
%about existence of $\rho$ for a given $\tau$ if $G$ is a compact group or a direct limit
%of compact groups.

%\sm  

{\bf \punct Reconstruction of a representation of $G$.%
\label{ss:reconstruction}}
To avoid an appearance of exotic measures let us fix a class of topological groups
under considerations.  Recall that a topological space is  Polish
if it is homeomorphic to a complete metric space.  
 A topological group is {\it Polish}
if it is a Polish topological space.

Consider a sequence 
$$
G_1\subset G_2\subset \dots
$$
of Polish groups. We say that the group
$$
G:=\cup G_j=:\lim_{\longrightarrow} G_j
$$ 
is their {\it direct limit}. We equip $G$ with the topology of a direct limit, a set
$U\subset G$ is open if all intersections $U\cap G_j$ are open in $G_j$. Generally speaking
such topologies are not metrizable. If a subset $K\subset G$ is compact, then 
 $K\subset G_j$ for some $j$.
 
 \sm 
 
 {\sc Remark.} A Polish group $G$ is a direct limit of Polish groups,
 $G_j=G$. \hfill $\boxtimes$

\sm 

 Let $G=\lim G_j$ be  a direct limit of Polish groups.
 Denote by $\cM(G)$ the algebra of all compactly supported Borel  complex-valued measures on
 $G$, this algebra is a direct limit of algebras
 $$
 \cM(G)=\lim_{\longrightarrow} \cM(G_j),
 $$ 
 Notice that such measures are objects of classical measure theory on compact metrizable
 spaces (see, e.g, \cite{RS}, Chapter 4). We say that a sequence $\mu_j\in \cM(G)$ converges to $\mu$
 if there is a compact subset $L\subset G$ such that $\mu_j(G\setminus L)=0$
 for all $j$ and $\mu_j$ weakly converge to $\mu$ on $L$.

 Let $K\subset G$ be a compact subgroup, without loss of  generality
 we can assume $K\subset G_1$. Denote by $\cM(G\kos K)\subset \cM(G)$ the subalgebra of
 $K$-bi-invariant compactly supported measures. 
 
 \begin{proposition}
 	\label{pr:1}
 	Let $G$ be a direct limit of Polish groups.
 	If a $*$-representation $\tau$  of the algebra $\cM(G\kos K)$ in a Hilbert space $V$
 	can be represented as $\cRes_{G\kos K}(\rho)$ and $\cRes_{G\kos K}(\rho')$,
 	then $\rho\simeq \rho'$.  
 \end{proposition}

We formulate a stronger version of the statement including a way of a reconstruction
of $\rho$. Consider the homogeneous space
$K\setminus G$, denote by $z_0$ the base point of this space, i.e., the coset
$K\cdot 1$. 
 Define a kernel $\cL(z,u)$ on  $K\setminus G\times K\setminus G$ taking values
in bounded operators $V\to V$ by 
$$
 \cL(x,y)=\tau\bigl(\delta_K* h g^{-1}*\delta_K \bigr), \qquad \text{where $z_0 g=x$, $z_0h=y$}
$$
(the result does not depend on a choice of $g$, $h$). Notice that
the kernel is $G$-invariant,
$$
\cL(xr,yr)=\cL(x,y), \qquad \text{for $r\in G$.}
$$

Consider the space $\Delta(K\setminus G,\tau)$ of {\it finitely} supported functions  $K\setminus G\to V$.
For a vector $v\in V$ denote by $v\delta_z(x)$ the function, which equals
$v$ at the point $z$ and 0 at over points.
So, our space $\Delta(K\setminus G)$ consists of finite linear
combinations 
$$
\sum_{j=1}^n v_j \delta_{z_j} (x).
$$

The group $G$ acts in the space $\Delta(K\setminus G,\tau)$ by shifts
of the argument.
Define a sesquilinear linear form on $\Delta(K\setminus G,\tau)$
setting
%$$
%\la v \delta_{a}, w \delta_b\ra_{\Delta}:=\la \cL(a,b) v, w\ra_V,
%$$
%and
\begin{equation}
\Bigl\la \sum v_i \delta_{a_i}, \sum w_j \delta_{b_j}\Bigr\ra=
\sum_{i,j}   \bigl\la \cL(a_i,b_j) v_i, w_j\bigr\ra_V.
\label{eq:ses}
\end{equation}
If for all vectors $\sum v_i \delta_{a_i}$ we have 
\begin{equation}
\Bigl\la \sum v_i \delta_{a_i}, \sum v_i \delta_{a_i}\Bigr\ra=
\sum_{i,j}   \bigl\la \cL(a_i,a_j) v_i, v_j\bigr\ra_V\ge 0,
\label{eq:positive}
\end{equation}
then $\la\cdot,\cdot\ra$ is an inner product, we
get a structure of a pre-Hilbert space and  take the corresponding
Hilbert space $\ov{\Delta(K\setminus G,\tau)}$. Since the kernel is $G$-invariant,
shifts on $K\setminus G$ induce unitary operators in $\ov{\Delta(K\setminus G,\tau)}$.

\begin{proposition}
	\label{pr:2}
	Let $G$ be a direct limit of Polish groups.
	Let $\tau=\cRes_{G\kos K}(\rho)$. Then \eqref{eq:positive} holds and the representation of $G$
	in $\ov{\Delta(K\setminus G,\tau)}$ is equivalent to $\rho$. 
\end{proposition}

{\sc Proof.} Consider a unitary representation $\rho\in \URep(G)_K$ in a Hilbert
space $H$ and the corresponding representation $\tau$ of $\cM(G\kos K)$ in $V:=H^K$.
We consider the map $\Delta(G\kos K,\tau)\to H$ determined by
\begin{equation}
J:\, \sum_i v_i \delta_{z_0 g_i}\mapsto \rho(g_i^{-1}) v_i.
\label{eq:delta-delta}
\end{equation}
We have
\begin{multline*}
\bigl\la \rho(g^{-1})v,\,\rho(h^{-1}) w\bigr\ra_H=\bigl\la \rho(g^{-1})\,\rho(\delta_K) v,\,\rho(h^{-1})\, \rho(\delta_K) w\bigr\ra_H
=\\=
\bigl\la \rho(\delta_K) \,\rho(h)\, \rho(g^{-1})\,\,\rho(\delta_K) v, w\bigr\ra_H=
\bigl\la \rho(\delta_K *hg^{-1}*\delta_K) v,\, w\bigr\ra_H
=\\=\bigl\la \rho(\delta_K *hg^{-1}*\delta_K) v,\, w\bigr\ra_{H^K}.
\end{multline*}
Therefore
\begin{multline*}
\Bigl\la J\Bigl(\sum_i v_i \delta_{z_0 g_i}\Bigr), 
J\Bigl( w_j \delta_{z_0 g_j}\Bigr) \Bigr\ra_H=
\Bigl\la\sum_i \rho(g_i^{-1})v_i,\,\sum_j \rho(h_j^{-1}) w_j\Bigr\ra_H
=\\=
\sum_{i,j} \bigl\la \rho(\delta_K *h_jg_i^{-1}*\delta_K) v_i,\, w_j\bigr\ra_{H^K}
=
\\=
\sum_{i,j} \bigl\la \cL(z_0 g_i, z_0 h_j) v_i,\, w_j\bigr\ra_{H^K}.
\end{multline*}

So the map $J$ induces an inner product on the space $\Delta(K\setminus G,\tau)$,
and this inner product coincides with the sesquilinear form \eqref{eq:ses}. So
(\ref{eq:positive}) is positive and the map 
\eqref{eq:delta-delta} determines a
unitary operator $\ov{\Delta(K\setminus G,\tau)}\to H$.
\hfill $\square$

\sm 

{\sc Remark.}
If $\tau$ is an arbitrary $*$-representation of $\cM(G\kos K)$, then we can
repeat the definition of the space $\Delta(K\setminus G,\tau)$, but 
 positivity \eqref{eq:positive} of an operator-valued reproducing kernel $\cL(x,y)$ 
can be a heavy problem.
\hfill $\boxtimes$

\sm 

{\bf \punct  The case of compact groups.%
\label{ss:compact}}
%Let $L$ be a  compact group with countable base, let $\lambda$ be a probabilistic Haar
%measure on $L$. Denote by $\cC(G)$ the convolution
%algebra of continuous functions on $L$, equipped  with the  uniform topology.
%The involution is given by $f(g)\mapsto \ov{f(g^{-1})}$. 
%For a closed  subgroup $K\subset L$ we denote by $\cC(L\kos K)$
%the subalgebra in $\cC(G)$ consisting of $K$-bi-invariant functions.

%This algebra is a dense subalgebra
%in $\cM(L)$, for any $f\in\cC(G)$ we assign the measure $f\lambda$.

\begin{proposition}
	\label{pr:3}
	For a compact group $L$ with a countable base any $*$-re\-pre\-sen\-ta\-tion $\tau$ of the algebra $\cC(L\kos K)$
	has the form $\cRes_{L\kos K}(\rho)$. 
\end{proposition}

\begin{corollary}
	The same statement holds for  the algebra
	$\cM(L\kos K)$.
\end{corollary}
 
{\sc Proof.}
 Consider the  set
$\wh L$ of all irreducible representations 
$\rho_\alpha $ of the group $L$ defined up to equivalence%
\footnote{This set is finite if the group $L$ is finite, otherwise it is countable.},
denote by $H_\alpha$ 
(finite-dimensional) spaces
of these representations, by $\Mat(H_\alpha)$ the algebras of all operators in $H_\alpha$.

 Consider the Fourier transform on $L$, see, e.g., \cite{Kir}, Subsect. 12.2.
Namely, for
$f\in \cC(L)$  the operator-valued function
$$
\cF f(\alpha):=\rho_\alpha(f)\in \Mat(H_\alpha),\qquad \text{where $\alpha\in \wh L$},
$$
is called the Fourier transform $\cF f$ of $f$. %So an element of Fourier transform 
%is an element of the direct product of algebras $\Mat(Y_j)$.
By the definition, the Fourier transform sends convolutions to point-wise products.
The $\cF$ is injective, 
 a convenient
descriptions of the Fourier-image of $\cC(G)$ and of the induced
convergence in the Fourier-image are unknown. In any case,
functions $\alpha\mapsto \|\rho_\alpha(f)\|$ are bounded;
 the image of the Fourier transform
 contains
 the set of all finitary functions (i.e., functions whose elements
are zeros for all but a finite number of $\alpha$). Moreover, finitary functions
are dense in the Fourier-image. The uniform
convergence of functions $f_j\to f$ implies convergences
$\rho_\alpha(f_j)\to \rho_\alpha(f)$.

Denote by $\chi_\beta$ the character of $\rho_\beta$.
Then 
$$
\cF \chi_\beta(\alpha)=
\rho_\alpha(\chi_\beta)=
\begin{cases}
(\dim H_\beta)^{-1}\cdot 1, \qquad&\text{if $\alpha=\beta$;}
\\
0\qquad & \text{otherwise}.
\end{cases} 
$$

Next, consider the subset $\wh L_K\subset \wh L$ consisting of representations
having non-zero $K$-fixed vectors. Define the  Fourier transform 
$\cF_K$ on $\cC(L\kos K)$ as  operator-valued function
$$
\cF_K f(\alpha):=\wt\rho_\alpha(f)\in \Mat(H_\alpha^K).
$$
These operators are left upper blocks of matrices \eqref{eq:000}.
Properties of $\cF$ can easily translated to the corresponding properties
of $\cF_K$. In particular, the Fourier-image $\cF_K (\cC(L\kos K))$
contains a copy of each algebra $\Mat(H^K_\alpha)$, it consist of functions 
supported by one point $\alpha$.
 Notice that elements 
 $$\zeta_\beta:=\dim(H_\beta)\cdot \chi_\beta*\delta_K\in \cC(G)$$
are commuting idempotents,
$$
\zeta_\beta*\zeta_\gamma=\begin{cases}
\zeta_\beta, \qquad& \text{if $\beta=\gamma$};\\
0\qquad& \text{otherwise}
\end{cases}.
$$
Their Fourier-images are 
$$
\cF_K\zeta_\beta(\alpha)=\begin{cases}
1,\qquad &\text{if $\alpha=\beta$};
\\
0 \qquad &\text{otherwise}.
\end{cases}
$$

Now consider a $*$-representation $\tau$ of the algebra $\cC(G\kos K)\simeq \cF_K (\cC(G\kos K))$
in a Hilbert space $V$. Then operators $\tau(\zeta_\alpha)$ are commuting
projectors, denote by $V_\alpha$ their images, these subspaces are pairwise
orthogonal.
Moreover, $V=\oplus_\alpha V_\alpha$.

 Indeed, let $v\in (\oplus_\alpha V_\alpha)^\bot$. Then  it is annihilated
by all subalgebras $\Mat(H_\alpha^K)$ in $\cF_K (\cC(G\kos K))$.
But their sum is dense in the Fourier-image, therefore it is annihilated by
the whole algebra $\cF_K( \cC(G\kos K))$. Thus, $v=0$.

Thus, it is sufficient to construct a desired extension for each summand $V_\alpha$.
So, without loss of generality we can assume $V=V_\alpha$. Then for all $\beta\ne \alpha$
operators $\tau(\zeta_\beta)$ are zero and $\tau$ is zero on each $\Mat(H_\beta^K)$.
So
$$
\cC(L\kos K)\bigl/\ker \tau\simeq \Mat(H^K_\alpha).
$$
Any representation of a matrix algebra is a direct sum of irreducible  tautological representations.
By definition, each summand arises from the representation $\rho_\alpha$ of $L$.
\hfill $\square$

\sm

{\bf \punct  The case of direct  limits of compact groups.%
\label{ss:limit}}

\begin{proposition}
	\label{pr:4}
	Let $L$ be a direct limit
	$L:=\lim\limits_{\longrightarrow} L_j$
	of compact groups, let all $L_j$ have  countable bases of topology.
	Then  any $*$-representation $\tau$ of the algebra $\cM(L\kos K)$
	has the form $\cRes_{L\kos K}(\rho)$.
\end{proposition}

{\sc Proof.} It is sufficient to check the positivity of \eqref{eq:positive}.
Fix an expression
$\sum v_i \delta_{a_i}$. Since the summation is finite, we have only finite
collection $\{a_j\}\subset K\setminus L$. Therefore the subset is contained
in some $K\setminus L_j$. But the prelimit group $L_j$ is compact,
and we can apply Proposition \ref{pr:3}.
\hfill $\square$  

\section{Characters of $\cH_\infty(q)$}

\COUNTERS

{\bf \punct An algebra of operators in a tensor power.%
\label{ss:algebra}}
Consider the Hilbert space $V$ as in Subsect. \ref{ss:result} equipped with the basis $v_j$, the and the copy $W$
of $V$. Consider the tensor product $V\otimes W$ and a unit vector
$$
\xi:=\sum_{j\ne 0} a_j^{1/2} v_j\otimes w_j,\qquad \text{where $a_j\ge 0$,\, $\sum a_j=1$.} 
$$
Consider the tensor power
$$
(V\otimes W)^{\otimes n}\simeq V^{\otimes n}\otimes W^{\otimes n}
$$
 and the vector $\Xi:=\xi^{\otimes n}\in  (V\otimes W)^n$.
 Denote elements of  the natural orthonormal basis in our space by
 $$
 \eta\begin{bmatrix}
 i_1 &\dots& i_n\\ j_1&\dots&j_n
 \end{bmatrix}:=(v_{i_1}\otimes v_{j_1})\otimes \dots \otimes (v_{i_n}\otimes v_{j_n}),
 $$
 in this notation
 $$
 \Xi=\sum_{i_1, \dots, i_n\in \Z\setminus 0}\, \prod_{k=1}^n\, a_{i_k}^{1/2} \cdot \eta\begin{bmatrix}
 i_1 &\dots& i_n\\ i_1&\dots&i_n
 \end{bmatrix}.
 $$

We say that an operator in $(V\otimes W)^{\otimes n}$ is {\it $V$-diagonal}
if it has  the form 
$$
D^V(\Phi)\, \eta\begin{bmatrix}
i_1 &\dots& i_n\\ j_1&\dots&j_n
\end{bmatrix}=\Phi(i_1,\dots,i_n) \,\eta\begin{bmatrix}
i_1 &\dots& i_n\\ j_1&\dots&j_n
\end{bmatrix},
$$   
where $\Phi(i_1,\dots,i_n)$ is a bounded function $(\Z\setminus0)^n\to \C$.
    For  $\sigma\in S_n$ denote by $T^V(\sigma)$ the permutation of factors 
    $V$ in $(V\oplus W)^n$ 
    corresponding to $\sigma$,
$$
T(\sigma)\eta\begin{bmatrix}
i_1 &\dots& i_n\\ j_1&\dots&j_n
\end{bmatrix}=\eta\begin{bmatrix}
i_{\sigma^{-1}(1)} &\dots& i_{\sigma^{-1}(n)}\\ j_1&\dots&j_n
\end{bmatrix}.
$$    
    
     Denote by $\cA^V$ the algebra of operators in $(V\otimes W)^{\otimes n}$
    generated by permutations $T^V(\sigma)$ and  operators $D^V(\Phi)$.
    Any element of this algebra can be represented as a linear combination
    of the form
    $$
    \sum_{\sigma\in S_n} T^V(\sigma) D^V_\sigma(\Phi) \quad \text{or}\quad 
     \sum_{\sigma\in S_n}  D^V_\sigma(\wt\Phi) T^V(\sigma). 
    $$ 
In the same way we define an algebra $\cA^W$, it consists
of similar operators acting on the factors $W$. These two algebras commute.
The map
$$
\sum_{\sigma\in S_n} T^V(\sigma) D^V_\sigma \,\mapsto\, \sum_{\sigma\in S_n} T^W(\sigma) D^W_\sigma.
$$

We define a linear anti-automorphism ('transposition') in $\cA^V$ by 
$$
\Bigl(\sum_{\sigma\in S_n} T^V(\sigma) D^V_\sigma(\Phi)\Bigr)^t
=
\sum_{\sigma\in S_n} D^V_\sigma (\Phi) T^V(\sigma^{-1}). 
$$

\begin{proposition}
	\label{pr:ABBA}
	{\rm a)} For any elements $A^V$, $B^V\in \cA^V$ we have
	\begin{equation}
\bigl	\la A^V B^V\Xi,\Xi\bigl\ra=	\la B^V A^V\Xi,\Xi\bigr\ra.
	\end{equation}
	
	{\rm b)} For any element $B^W\in \cA^W$ we have
	\begin{equation}
	 B^W \Xi =  (B^V)^t \Xi.
	\end{equation}

	{\rm c)} For any elements $A^V\in \cA^V$, $B^W\in\cA^W$,
	\begin{equation}
\bigl	\la A^V B^W \Xi,\Xi\bigl\ra= \bigl 	\la (B^V)^t A^V  \Xi,\Xi\bigl\ra.
	\end{equation}
\end{proposition}

So we come to the case of Hilbert algebras discussed in Subsect. \ref{ss:double}.
The function 
$$
\chi(A^V)=\la A^V\Xi,\,\Xi\ra
$$
is a trace on the algebra
 $\cA^V$ (an explicit formula is \eqref{eq:TPhi}). The $\cA^V$-cyclic span
 of $\Xi$ is identified with $\ov{\cA^V}_\chi$, the map  $\cA^V\to \ov{\cA^V}_\chi$
 is $A^V\to A^V \Xi$. The right action of our algebra on $\ov{\cA^V}_\chi$
 is given by $X\mapsto X (A^W)^t$. 
 
 \sm  
 
 The construction of Subsect. \ref{ss:result}
 gives an embedding 
 $$\cH_{\infty}(q)\to \cA^V$$
 and an action of $\cH_{\infty}(q)\otimes \cH_{\infty}(q)$ on the $\ov{\cA^V}_\chi$.
 In Subsect. \ref{ss:zeta} we will show that this embedding induces the Vershik--Kerov 
 trace.
 
 \sm

{\bf \punct Proof of Proposition \ref{pr:ABBA}.%
\label{ss:algebra-proof}}
We use notation $I:=(i_1,\dots,i_n)$, $J:=(j_1,\dots,j_n)$.
For $\sigma\in S_n$ denote by $\Omega(\sigma)$ the set of all $I$ invariant with respect
to $\sigma$. In other words we decompose $\sigma$ into a product of independent cycles,
and map $k\mapsto i_k$ is constant on cycles.

\begin{lemma}
	\label{l:TPhi}
	\begin{multline}
	\bigl\la   T^V(\sigma) D^V(\Phi)\,\Xi,\Xi\bigr\ra=
	\sum_{I=(i_1,\dots,i_n)\in \Omega(\sigma)} \Bigl( \prod_{k=1}^n a_{i_k}\cdot \Phi(I)\Bigr)
	=\\=\Bigl\la D^V(\Phi) T^V(\sigma)\, \Xi, \Xi\Bigr\ra.
	\label{eq:TPhi}
	\end{multline}
\end{lemma}

 {\sc Proof.}
We have
\begin{multline}
\Bigl\la  T^V(\sigma) D^V(\Phi) \Xi, \Xi\Bigr\ra=
\sum_I \sum_J\, \boxed{\Phi(i_{\sigma^{-1}(1)}, \dots, i_{\sigma^{-1}(n)}) }\,
\prod_{k=1}^n a_{i_{\sigma^{-1}(k)}}^{1/2} \prod_{k=1}^n a_{j_k}^{1/2}
\times \\ \times
\biggl\la \eta\begin{bmatrix}
i_{\sigma^{-1}(1)}& \dots& i_{\sigma^{-1}(n)}\\
i_1&\dots&i_n
\end{bmatrix},
\eta\begin{bmatrix}
j_1&\dots &j_n\\j_1&\dots &j_n
\end{bmatrix}
\biggr\ra.
\label{eq:Phi-boxed}
\end{multline}
 A summand can be non-zero only if two basis elements in the brackets $\la\cdot,\cdot\ra$
 coincide. So $i_k=j_k$, $i_{\sigma^{-1}k} =j_k$. Therefore $\sum_I\sum_J$
 transforms to $\sum_{I\in \Omega(\sigma)}$.
 Since $I\in \Omega(\sigma)$, we can replace $\Phi(\cdot)$ by $\Phi(I)$.
 Since $J=I$, the product $\prod_k(\dots)\prod_k(\dots)$ replaces by
 $\prod_k a_{i_k}$.
 
 The expression for $\bigl\la D^V(\Phi) T^V(\sigma)\, \Xi, \Xi\bigr\ra$
 is similar to \eqref{eq:Phi-boxed}, we only must replace the boxed
 $\Phi(\dots)$ by $\Phi(i_1,\dots,i_n)$.
 \hfill $\square$

 \begin{lemma}
 	\label{l:PhiTheta}
 	\begin{equation}
 	\Bigl \la \bigl(T^V(\sigma) D^V(\Phi)\bigr)\,D^V(\Theta)\Xi,\, \Xi\Bigr\ra=
 	\Bigl\la D^V(\Theta)\,  \bigl( T^V(\sigma) D^V(\Phi)\bigr)\Xi,\, \Xi\Bigr\ra.
 	\label{eq:PhiTheta}
 	\end{equation}
 \end{lemma}

{\sc Proof.} 
To obtain the left hand side of \eqref{eq:PhiTheta}, we must replace the boxed $\Phi(\dots)$
in the calculation \eqref{eq:Phi-boxed} by
$$\Phi(i_{\sigma^{-1}(1)}, \dots, i_{\sigma^{-1}(n)})\,\Theta(i_{\sigma^{-1}(1)}, \dots, i_{\sigma^{-1}(n)}),$$
in the right hand side by
$$\Theta(  i_1, \dots, i_n)\,\Phi(i_{\sigma^{-1}(1)}, \dots, i_{\sigma^{-1}(n)} ). $$
Since $I\in \Omega(\sigma)$, we will get the same result in both sides.
\hfill $\square$

\begin{lemma}
	\label{l:lambda-sigma}
	\begin{equation}
	\Bigl\la T^V(\lambda^{-1})\, T^V(\sigma) D^V(\Phi)\, T^V(\lambda)\, \Xi,\,\Xi \Bigr\ra=
	\Bigl\la  T^V(\sigma) D^V(\Phi)\, \Xi,\,\Xi \Bigr\ra.
	\end{equation}
\end{lemma}

{\sc Proof.} The left-hand side equals to
\begin{multline*}
\sum_I\sum_J \Phi(i_{\sigma^{-1}\lambda(1)},\dots,i_{\sigma^{-1}\lambda(n) }) 
	\prod_{k=1}^n a_{i_{\lambda^{-1}\sigma^{-1} \lambda(k)}} \prod_{k=1}^n a_{j_k}
	\times \\ \times
	\times \\ \times
	\biggl\la \eta\begin{bmatrix}
	i_{\lambda^{-1}\sigma^{-1}\lambda(1)}& \dots& i_{\lambda^{-1}\sigma^{-1}\lambda(n)}\\
	i_1&\dots&i_n
	\end{bmatrix},
	\eta\begin{bmatrix}
	j_1&\dots &j_n\\j_1&\dots &j_n
	\end{bmatrix}
	\biggr\ra.
\end{multline*}
As in proof of Lemma \eqref{eq:TPhi}, we come to
$$
\sum_{I\in \Omega(\lambda^{-1}\sigma^{-1}\lambda)}
\Phi(i_{\sigma^{-1} \lambda(1)},\dots,i_{\sigma^{-1}\lambda(n) })\prod_{k=1}^{n} a_{i_k}.
$$
(we applied Lemma \ref{l:TPhi}).
Next, we pass to a new index of summation 
$$I':=(i_{\lambda(1)}, \dots,i_{\lambda(n)}) ,\qquad I'\in \Omega(\sigma).$$ 
Keeping in mind Lemma \ref{l:TPhi} we get the desired expression in the right hand side.
\hfill $\square$

\sm 

The statement a) of Lemma \ref{pr:ABBA} follows from
Lemmas \ref{l:PhiTheta} and \ref{l:lambda-sigma}.

\sm 

{\sc Proof of the statement b) of Proposition \ref{pr:ABBA}.}
%We must verify the identity
%$$
%\Bigl\la T^V(\sigma)D^V(\Phi)\, T^W(\lambda) D^W(\Psi)\,\Xi,\Xi \Bigr\ra=
%\Bigl\la D^V(\Psi) T^V(\lambda^{-1})\, T^V(\sigma) D^V(\Phi)\,\Xi,\,\Xi\Bigr\ra.
%$$
We must verify the identity
$$
 T^W(\lambda) D^W(\Psi)\,\Xi=
 D^V(\Psi) T^V(\lambda^{-1})\,\Xi.
$$
In the left hand-side we have
\begin{equation*}
 \sum_I  \Psi(i_{\lambda^{-1}(1)},\dots,i_{\lambda^{-1}(1)})  \prod a^{-1/2}_{j_{\lambda^{-1} k}}
\, \eta
\begin{bmatrix}
i_1&\dots&i_n\\
i_{\lambda^{-1}(1)}& \dots& i_{\lambda^{-1}(n)}
\end{bmatrix}.
\end{equation*}

%The left hand side equals to
%\begin{multline*}
%\sum_I \sum_J  \Phi(i_{\sigma^{-1}(1)},\dots, i_{\sigma^{-1}(n)}) \Psi(i_{\lambda^{-1}(1)},\dots,i_{\lambda^{-1}(1)}) \prod_k a^{-1/2}_{i_{\sigma^{-1} k}} \prod a^{-1/2}_{j_{\lambda^{-1} k}}
%\times\\\times
%\biggl\la \eta
%\begin{bmatrix}
%i_{\sigma^{-1}(1)}&\dots&i_{\sigma^{-1}(n)}\\
%i_{\lambda^{-1}(1)}& \dots& i_{\lambda^{-1}(n)}
%\end{bmatrix},
%\eta\begin{bmatrix}
%j_1,\dots,j_n\\j_1,\dots,j_n
%\end{bmatrix}
%\biggr\ra 
%=\\=
%\sum\limits_{(i_1,\dots,i_n):\,i_{\sigma^{-1}(k)}=i_{\lambda^{-1}(k)}}
%\Phi(i_{\sigma^{-1}(1)},\dots, i_{\sigma^{-1}(n)}) \Psi(i_{\lambda^{-1}(1)}, \dots, i_{\lambda^{-1}(n)}) \prod_k a_{i_{\sigma^{-1}(k)}}.
%\end{multline*}
%The set of summation is $I\in \Omega(\lambda \sigma^{-1})=\Omega(\sigma\lambda^{-1})$,
%arguments of $\Phi$ and $\Psi$ coincide.
 We pass to a new index of summation
$$
I':=(i_{\lambda^{-1}(1)}, \dots, i_{\lambda^{-1}(n)}),%\qquad
%I'\in \Omega(\lambda^{-1}\sigma)
$$
and get the right hand side.
\hfill $\square$

\sm

The statement c) of Proposition \ref{pr:ABBA} follows from b).

\sm 

{\bf \punct The matrix element of $\zeta_{m}$.%
\label{ss:zeta}}

\begin{lemma}
	\label{l:zeta}
	Let $\cH_n$ acts in $(V\otimes W)^{\otimes n}$ as in Subsect. {\rm\ref{ss:limit}}.
	Then 
	$$
	\bigl\la \zeta_{m}\Xi,\Xi\ra=\chi^{\alpha,\beta,0}(\zeta).
	$$
\end{lemma}

\begin{lemma}
	\label{l:RQD}
	Decompose the $R$-matrix \eqref{eq:R} as
	$$R:=Q+D,$$
	where
	\begin{align*}Q:&=-\sqrt q \sum_{i\ne j,i\ne0, j\ne0} e_{ji}\otimes f_{ij};
	\\
	D:&=-\sum_{i<0} e_{ii}\otimes f_{ii} +q \sum_{i>0} e_{ii}\otimes f_{ii}
	 +(q-1) \sum_{i< j,i\ne0, j\ne0}
	e_{ij}\otimes f_{ij}.
	\end{align*}
	Then
	$$
	R_{(m-1)m}\dots R_{23}R_{12}\,\Xi=D_{(m-1)m}\dots D_{23}D_{12}\,\Xi.
	$$	
\end{lemma}

{\sc Proof of Lemma \ref{l:RQD}.}
Denote 
$$a_{i}:=\begin{cases}\beta_{-i},\quad \text{for $i<0$};\\
\alpha_i,\quad \text{for $i>0$.}
\end{cases}
$$
The operators $D_{j(j+1)}$ are $V$-diagonal. Next,
$$
Q_{12}\Xi =-\sqrt q \sum_{I=(i_1,\dots,i_n): i_1\ne i_2}
\prod_k a_{i_k}^{1/2} \eta\begin{bmatrix}
i_2& i_1& i_3&\dots &i_n\\ 
i_1& i_2& i_3&\dots &i_n
\end{bmatrix}
$$
The operators $R_{23}$, $R_{34}$, \dots
can not change $i_2$ in the first column of $\eta[\cdot]$,
recall also that such operators do not act on the second row of $\eta[\cdot]$.
Therefore all terms of $R_{(m-1)m}\dots R_{23}Q_{12}\,\Xi$
have form
$$
c_s \eta \begin{bmatrix}
\alpha_1^s&\dots&\alpha_n^s\\
\beta_1^s&\dots&\beta_n^s
\end{bmatrix},\qquad \text{where $\alpha_1^s\ne \beta_1^s$.}
$$
 Inner product of such a term with $\Xi$ is 0, i.e.,
$$
\bigl\la R_{(m-1)m}\dots R_{23}Q_{12}\,\Xi,\Xi \bigr\ra=0,
$$
hence
$$
\bigl\la R_{(m-1)m}\dots R_{23}R_{12}\,\Xi,\Xi \bigr\ra=
\bigl\la R_{(m-1)m}\dots R_{23}D_{12}\,\Xi,\Xi \bigr\ra.
$$
We repeat the same argument for $R_{23}$, $R_{34}$, etc.
and get the desired statement.
\hfill $\square$

\sm

The operator $D_{(m-1)m}\dots D_{23}D_{12}$ is $V$-diagonal. Denote by $\delta[\dots]$ its eigenvalues
\begin{multline*}
D_{(m-1)m}\dots D_{23}D_{12}\,\eta\begin{bmatrix}
i_1&\dots&i_m&\dots  &i_n\\j_1&\dots&j_m&\dots &j_n
\end{bmatrix}
=\\= \delta(i_1,\dots,i_m)\,
\eta\begin{bmatrix}
i_1&\dots&i_m&\dots  &i_n\\j_1&\dots&j_m&\dots  &j_n
\end{bmatrix}.
\end{multline*}

\begin{lemma}
	\begin{equation}
	\la \zeta_m \Xi,\Xi\ra= \sum_{i_1,\dots,i_m} \Bigl(\delta(i_1,\dots,i_m) \prod_{k=1}^m a_{i_k}\Bigr) .
	\label{eq:zeta}
		\end{equation}
\end{lemma}

{\sc Proof.} A straightforward summation gives
$$
\sum_{i_1,\dots,i_n} \Bigl(\delta(i_1,\dots,i_m) \prod_{k=1}^n a_{i_k}\Bigr) 
$$
We transform this expression as
\begin{multline*}
\sum_{i_1,\dots,i_m} \sum_{i_{m+1},\dots,i_n}
\Bigl(\delta(i_1,\dots,i_m) \prod_{k=1}^m a_{i_k} \prod_{k=m+1}^n a_{i_k}\Bigr)
=\\=
 \sum_{i_1,\dots,i_m} \Bigl(\delta(i_1,\dots,i_m) \prod_{k=1}^m a_{i_k}\Bigr) 
 \times \Bigl(\sum_i a_i\Bigr)^{m-n}.
\end{multline*}
 But under our conditions $\sum a_i =\sum \beta_p+\sum \alpha_q=1$.
 \hfill $\square$
 
 \sm 

Next, $\delta(i_1,\dots,i_n)$ is non-zero iff
$i_1\le \dots \le i_n$. The next statement also is obvious:

\begin{lemma} Let a tuple $I:i_1\le \dots \le i_n$
consists of entries $(-\iota_u),\dots, (-\iota_1)<0$ with nonzero multiplicities 
$\mu_1$, \dots, $\mu_u$ and entries $\upsilon_1$, \dots, $\upsilon_v$
with nonzero multiplicities $\nu_1$, \dots, $\nu_v$. Then 	
\begin{equation}
\delta(i_1,\dots,i_n)=
(-1)^{\sum_k (\mu_k-1)} q^{\sum_l(\nu_l-1)} (q-1)^{u+v-1}.
\label{eq:delta}
\end{equation}
\end{lemma}

\begin{corollary}
	\begin{multline}
	\la \zeta_m \Xi,\Xi\ra= \frac1{q-1}
	\sum_{\text{$\begin{matrix}\phi_1\ge 0,\phi_2\ge 0, \dots,\\
	\psi_1\ge 0,\psi_2 \ge0,\dots:\\ 
	\sum \phi_i+\sum \psi_j=m\end{matrix}$}}\,
\prod_{i=1}^\infty \Bigl((-\beta_i)^{\phi_i} (1-q)^{\epsilon(\phi_i)}\Bigr)
\times \\\times 
\prod_{j=1}^{\infty} \Bigl( (q\alpha_j)^{\psi_l}(1-q^{-1})^{\epsilon(\psi_i)} \Bigr),
\label{eq:phiphi}
	\end{multline}
	where
	$$
	\epsilon(\theta):=\begin{cases}
	0, \qquad \text{if $\theta=0$;}
	\\
	1,\qquad \text{if $\theta>0$.}
	\end{cases}
	$$
\end{corollary}

{\sc Proof.} We transform the expression \eqref{eq:delta} for $\delta(\dots)$ as
$$
\frac{1}{q-1}\prod_k \Bigl( (-1)^{\mu_k} (1-q)\Bigr)\cdot \prod_l  
\Bigl( q^{\nu_l} (1-q^{-1}) \Bigr)
.
$$
 A straightforward summation  in formula \eqref{eq:zeta} 
gives an expression
\begin{multline*}
\la \zeta_m \Xi,\Xi\ra=\frac{1}{q-1}\sum_{u,v:u+v\ge 1}\,\,
\sum_{\iota_1<\dots <\iota_u, \upsilon_1<\dots<\upsilon_v}\,\,
\sum_{\mu_1,\dots,\mu_u,\nu_1,\dots,\nu_v:\sum \mu_k+\sum \nu_l=m}
\times \\\times 
\prod_{k=1}^u \Bigl((-\beta_{\iota_k})^{\mu_k} (1-q)\Bigr)
\prod_{k=1}^{v} \Bigl( (q\alpha_{\upsilon_l})^{\nu_l}(1-q^{-1}) \Bigr).
\end{multline*}
We change a parametrization of the set of summation assuming
\begin{align*}
\phi_{\iota_k}&=\mu_k \text{and $\phi_i=0$ for all other $i$;}
\\
\psi_{\upsilon_l}&=\nu_l  \text{and $\psi_j=0$ for all other $l$,}
\end{align*}
and get \eqref{eq:phiphi}, all new factors in the products in
\eqref{eq:phiphi} are 1.
\hfill $\square$

\begin{lemma}The generating function 
	$$G(z):=1+(q-1)\Bigl(z +\sum_{m\ge 2} \la \zeta_m \Xi,\Xi\ra z^m\Bigr)$$
	is equal to
	\begin{equation}
G(z):=	\prod_{i=1}^{\infty} \frac{1+\beta_i qz}{1+\beta_i z}
	\,\prod_{j=1}^{\infty}\frac{1-\alpha_j z}{1-\alpha_j q z}.
	\label{eq:generating}
		\end{equation}
\end{lemma}

{\sc Proof.} We represent the term  $(q-1)z$	
as 
$$
(q-1)z=(q-1)\Bigl(\sum_j \beta_j+\sum_k \alpha_k\Bigr)z=\sum_j (-\beta_jz)(1-q)+\sum_k (\alpha_k qz)(1-q^{-1}).
$$
Applying \eqref{eq:phiphi} we come to the following expression for the generating function:
	\begin{multline*}
G(z)=\sum_{\phi_1\ge 0,\phi_2\ge 0, \dots,
	\psi_1\ge 0,\psi_2\ge0,\dots }
\prod_{i=1}^\infty \Bigl((-\beta_i z)^{\phi_i} (1-q)^{\epsilon(\phi_i)}\Bigr)
\times \\\times 
\prod_{j=1}^{\infty} \Bigl( (q\alpha_j z)^{\psi_j}(1-q^{-1})^{\epsilon(\psi_j)} \Bigr).
\end{multline*}
It decomposes into a product
\begin{multline*}
\prod_{i=1}^\infty \Bigl(\sum_{\phi_i}(-\beta_i z)^{\phi_i} (1-q)^{\epsilon(\phi_i)}\Bigr)
\cdot 
\prod_{j=1}^{\infty} \Bigl(\sum_{\psi_j} (q\alpha_j z)^{\psi_j}(1-q^{-1})^{\epsilon(\psi_j)} \Bigr)
=\\=
\prod_{i=1}^\infty \Bigl(1+(1-q)\bigl((1+\beta_i z)^{-1}-1 \bigr) \Bigr)
\,
\prod_{j=1}^{\infty} \Bigl(1+(1-q^{-1})\bigl( (1-q\alpha_j z)^{-1}-1\bigr) \Bigr).
\end{multline*}
This equals to \eqref{eq:generating}.
\hfill $\square$

\sm 

{\sc Proof of Lemma \ref{l:zeta}.}	We must verify the identity
$$
G(z)=1+(q-1)z+\sum_{m\ge 2} z^m \chi^{\alpha,\beta,0}(\zeta_m),
$$
see \eqref{eq:vk}. We have
\begin{multline*}
G(z)=\exp\{ \ln G(z)   \}=\\=
\exp\Bigl\{\sum_i\bigl(\ln(1+\beta_i q z)-\ln(1+\beta_iz)\bigr)+
\sum_j \bigl(\ln(1-\alpha_j z)-\ln(1-\alpha_j q z)\bigr)
 \Bigr\}=\\=
 \exp\Bigl\{\sum_i\sum_k\frac{(-1)^{k-1}}k  \beta_i^k(q^k-1)z^k
+\sum_j \sum_k\frac 1k\, \alpha_j^k (q^k-1) z^k
  \Bigr \}
  =\\=
\prod_{k=1}^{\infty}  \exp\Bigl\{\bigl(\sum_j \alpha_j^k+(-1)^{k-1} 
  \sum_j\beta_i^k \bigr)\cdot \frac 1k\,  (q^k-1) z^k  \Bigr \}.
\end{multline*}
It remains to decompose exponentials.
\hfill $\square$

\sm 

{\bf \punct Proof of Theorem \ref{th:}.}
Obviously, $R_{i(i+1)}^{\mathrm{left}}$ are contained in the algebra
$\cA^V$ introduced in Subsect. \ref{ss:algebra}. Therefore
operators of the representation $\cH_n(q)$ are contained in $\cA^V$.
By Proposition \ref{pr:ABBA}.a, the matrix element $\la A\Xi,\Xi\ra$ is a trace
on each $\cH_n$.  Therefore it is a trace on $\cH_\infty(q)$.
The properties \eqref{eq:multi}--\eqref{eq:shift} are obvious.
After these remarks the statement a) of Theorem \ref{th:} follows from Lemma \ref{l:zeta}
and the statement b) from Proposition \ref{pr:ABBA}.b.

 \tt
\noindent
Yury Neretin\\
 Math. Dept., University of Vienna \\
\&Institute for Theoretical and Experimental Physics (Moscow); \\
\&MechMath Dept., Moscow State University;\\
\&Institute for Information Transmission Problems;\\
yurii.neretin@math.univie.ac.at
\\
URL: http://mat.univie.ac.at/$\sim$neretin/

\end{document}